\documentclass[10pt]{article}
\usepackage{amsmath}
\usepackage{amsfonts}
\usepackage{amssymb}
\usepackage{comment}
\usepackage{graphicx,color}
\long\def\proof#1{\removelastskip\vskip\baselineskip\relax\noindent{\it
Proof\if!#1!\else\ \ignorespaces#1\fi.\ }\ignorespaces}

\newcommand{\ov}{\overline}
\newcommand{\Q}{{\mathbb Q}}
\newcommand{\Z}{{\mathbb Z}}

\newcommand{\z}{\zeta}

\newcommand{\ga}{\gamma}

\newcommand{\Proof}{{\it Proof. \/}}
\newcommand{\squareforqed}{\hbox{\rlap{$\sqcap$}$\sqcup$}}
\newcommand{\qed}{\ifmmode\squareforqed\else{\unskip\nobreak\hfil
\penalty50\hskip1em\null\nobreak\hfil\squareforqed
\parfillskip=0pt\finalhyphendemerits=0\endgraf}\fi}

\newcommand{\fp}{\qed\removelastskip\vskip\baselineskip\relax}

\newtheorem{theorem}{Theorem}[section]
\newtheorem{corollary}[theorem]{Corollary}
\newtheorem{proposition}[theorem]{Proposition}
\newtheorem{lemma}[theorem]{Lemma}
\newtheorem{definition}[theorem]{Definition}

\newcommand{\litem}{\par\noindent\dimen0=\parindent%
    \advance\dimen0 by-4pt
               \hangindent=\dimen0\ltextindent}

\newcommand{\ltextindent}[1]{\hbox to \hangindent{#1\hss}\ignorespaces}
\newcommand{\ltextjndent}[1]{\hbox to \hangindent{#1\hss}\ignorespaces\kern-1ex}

\renewcommand{\pmod}[1]{\allowbreak\ ({\rm{mod}}\,\,#1)}

\begin{document}
\pagestyle{plain}

\title{Elementary Continued Fractions for Linear Combinations of Zeta and $L$ Values}
\author{Henri Cohen,\\ henri.cohen2@free.fr}

\maketitle
\begin{abstract}
  We show how to obtain infinitely many continued fractions for certain
  $\Z$-linear combinations of zeta and $L$-values. The methods are completely
  elementary.
\end{abstract}

\smallskip

\section{Introduction and Main Goal}

In \cite{Ram}, the authors give five continued fractions for certain
$\Z$-linear combinations of zeta values, only obtained and checked numerically,
as well as other linear combinations involving powers of $\pi$, Catalan's
constant, etc... The purpose of the present paper is to show that these
continued fractions are completely elementary. In fact, we explain three
general methods for constructing them. In particular we give seven infinite
families of such continued fractions. We also discuss analogous results
involving $L$-values for Dirichlet characters of conductor $3$ and $4$.

We use the following notation for a continued fraction $S$, which may
differ from notation used in other papers in the literature:
$$S=a(0)+b(0)/(a(1)+b(1)/(a(2)+b(2)/(a(3)+\cdots)))\;,$$
and we denote as usual by $p(n)/q(n)$ the $n$th convergent, so that
$p(0)/q(0)=a(0)/1$, $p(1)/q(1)=(a(0)a(1)+b(0))/a(1)$, etc...

When $a(n)$ and $b(n)$ are polynomials $A(n)$ and $B(n)$ for $n\ge1$, we will
write the continued fraction as $S=[[a(0),A(n)],[b(0),B(n)]]$.
For instance, the continued fraction
$$\z(2)=2/(1+1/(3+16/(5+81/(7+256/(9+\cdots)))))$$
will simply be written as $\z(2)=[[0,2n-1],[2,n^4]]$.

\smallskip

We recall the following trivial result due to Euler:
\begin{lemma} Let $f(n)$ be a nonzero arithmetic function, and set by
  convention $f(0)=0$. When the left-hand side converges, we have
  $$\sum_{n\ge1}\dfrac{z^n}{f(n)}=[[0,f(n)+zf(n-1)],[z,-zf(n)^2]]\;,$$
  and in addition the $N$th partial sum of the series is equal to
  the $N$th convergent $p(N)/q(N)$ of the continued fraction.
\end{lemma}

As a trivial application, for $k\ge2$, we have the trivial continued fraction
$$\z(k)=[[0,n^k+(n-1)^k],[1,-n^{2k}]]\;.$$

The second trivial lemma that we will use is the following:
\begin{lemma} Let $(a(n),b(n))$ define a continued fraction with convergents
  $(p(n),q(n))$, and let $r(n)$ be an arbitrary nonzero arithmetic function
  with $r(0)=1$. Then if we set $a'(n)=r(n)a(n)$ and $b'(n)=r(n)r(n+1)b(n)$,
  the corresponding convergents $(p'(n),q'(n))$ are given by
  $(p'(n),q'(n))=r!(n)(p(n),q(n))$ with evident notation, and in particular
  $p'(n)/q'(n)=p(n)/q(n)$.\end{lemma}

Thanks to the first lemma above, we can thus transform any series into
a continued fraction, not really interesting. For instance, assume that I
want a CF for $\z(2)+\z(3)$: we have $\z(2)+\z(3)=\sum_{n\ge1}(n+1)/n^3$,
so applying the first lemma to $f(n)=n^3/(n+1)$ and $z=1$, we get
$$\z(2)+\z(3)=[[0,n^3/(n+1)+(n-1)^3/n],[1,-n^6/(n+1)^2]]\;,$$
and applying the second lemma to $r(0)=1$, $r(n)=n^2+n$ for $n\ge1$, we get
$$\z(2)+\z(3)=[[0,2n^4-2n^3+2n-1],[2,-(n^8+2n^7)]]\;.$$
As mentioned, not very interesting, in particular because continued
fractions involving $\z(3)$ should have $a(n)$ a polynomial of degree at most
$3$, and $b(n)$ at most $6$.

\smallskip

We can now more precisely state our goal, much wider than the simple proofs
of the continued fractions given in \cite{Ram}. First, we set the following
definition:
\begin{definition} Let $d\ge0$ be an integer. A \emph{rational period} of
  degree $k$ is the sum of a convergent series of the form
  $\sum_{n\ge1}\chi(n)f(n)$, where $\chi(n)$ is a periodic arithmetic function
  taking rational values, and $f\in\Q(x)$ is a rational function with
  rational coefficients, whose denominator is of degree $k$.\end{definition}

Two remarks concerning this definition: first, it is \emph{not} compatible
with the definition of periods as given in \cite{Kon-Zag}. Second, one
could ask the coefficients to be algebraic instead of rational, but this
leads to a theory which is too general.

Examples: $\log(2)$, $\pi$, and $L(\chi,1)$ for a nontrivial Dirichlet
character all have degree $1$, more generally $\pi^k$ and $L(\chi,k)$
are rational periods of degree $k$. Note also that a $\Q$-linear combination
of rational periods of degree $\le k$ is again a rational period of
degree $\le k$.

Thus, our goal will be as follows: find continued fractions $(a(n),b(n))$
for rational periods degree $k$ where for $n$ sufficiently large $a(n)$ is a
polynomial of degree at most $k$ and $b(n)$ of degree at most $2k$,
which we abbreviate by saying that it has \emph{bidegree} at most $(k,2k)$.

The prototypical ``trivial'' example is $S=\sum_{n\ge1}1/P(n)$ with
$P$ a polynomial of degree $k\ge2$, and $S=[[0,P(n)+P(n-1)],[1,-P(n)^2]]$
is a continued fraction of required bidegree $(k,2k)$, by Euler's lemma
above, and similarly for $S=\sum_{n\ge1}(-1)^{n-1}/P(n)$.
On the contrary, as mentioned above, such a continued fraction does not
seem to exist for $\sum_{n\ge1}(n+1)/n^3$.

\section{First Method: use of Polynomial Multipliers}

\begin{proposition} Fix an integer $k\ge2$, and let $P\in\Q[x]$ be a nonzero
  polynomial with rational coefficients such that $P(x)$ divides
  $x^kP(x+1)+(x-1)^kP(x-1)$, and set $R(x)=(x^kP(x+1)+(x-1)^kP(x-1))/P(x)$.
  \begin{enumerate}\item
  We have the continued fraction expansion
  $$S=\sum_{n\ge1}\dfrac{1}{n^kP(n)P(n+1)}=[[0,R(x)],[1/P(0)^2,-n^{2k}]]\;,$$
    which is a continued fraction of bidegree $(k,2k)$.
  \item If $P(x)$ and $P(x+1)$ are coprime polynomials and $d=\deg(P)\le k$,
    $S$ is a rational period of degree at most $k$.
  \end{enumerate}
\end{proposition}

\Proof (1). By Euler's lemma above we have
$$S=[[0,n^kP(n)P(n+1)+(n-1)^kP(n-1)P(n)],[1,-n^{2k}P(n)^2P(n+1)^2]]\;.$$
I claim that $P(1)\ne0$: indeed, if $P(1)=0$ we deduce that $1^kP(2)=R(1)P(1)$
so $P(2)=0$, and the recursion $(n-1)^kP(n)=R(n-1)P(n-1)-(n-2)^kP(n-2)$
for $n\ge3$ implies that $P(n)=0$ for all $n$, so $P$ has infinitely many
roots so is identically zero, contradiction. We can thus apply the second
lemma to $r(n)=1/P(n)^2$ for $n\ge1$ and $r(0)=1$ and
we obtain immediately $S=[[0,R(n)],[1/P(1)^2,-n^{2k}]]$, proving (1).

\smallskip

For (2), we first note that using the same proof as in (1) but using the
recursion backwards, we have $P(0)\ne0$. We can thus write a partial fraction
decomposition in the form
$$\dfrac{1}{x^kP(x)P(x+1)}=\sum_{1\le j\le k}\dfrac{c_j}{x^j}+\dfrac{N(x)}{P(x)P(x+1)}\;,$$
where the $c_j$ are constants and $\deg(N(x))\le 2d-1$, where $d=\deg(P)$.
Since by assumption $P(x)$ and $P(x+1)$ are coprime, by the extended Euclidean
algorithm, there exist polynomials $U$ and $V$ such that
$U(x)P(x)+V(x)P(x+1)=N(x)$, and $U$ and $V$ can be chose of degree less than
or equal to $d-1$. Thus $N(x)/(P(x)P(x+1))=U(x)/P(x+1)+V(x)/P(x)$, so
$S$ is a rational period of degree at most $d\le k$.\fp

\smallskip

{\bf Remarks}

\smallskip

\begin{enumerate}\item It is possible that the condition that $P(x)$ and
  $P(x+1)$ are coprime can be lifted.
\item One can prove that $P(x)$ satisfies the identity $P(x+1)=(-1)^dP(-x)$,
  or equivalently $P(1-x)=(-1)^dP(x)$. I am grateful to ``Ilya Bogdanov''
from the MathOverflow forum for the proof of this fact.\end{enumerate}

\medskip

Exactly the same proposition with an identical proof can be applied to
alternating sums:

\begin{proposition} Fix an integer $k\ge2$, and let $P\in\Q[x]$ be a nonzero
  polynomial with rational coefficients such that $P(x)$ divides
  $x^kP(x+1)-(x-1)^kP(x-1)$, and set $R(x)=(x^kP(x+1)-(x-1)^kP(x-1))/P(x)$.
  \begin{enumerate}\item
  We have the continued fraction expansion
  $$S=\sum_{n\ge1}\dfrac{(-1)^{n-1}}{n^kP(n)P(n+1)}=[[0,R(x)],[1/P(0)^2,n^{2k}]]\;,$$
    which is a continued fraction of bidegree $(k-1,2k)$.
  \item If $P(x)$ and $P(x+1)$ are coprime polynomials and $d=\deg(P)\le k$,
    $S$ is a rational period of degree at most $k$.
  \end{enumerate}
\end{proposition}

\smallskip

The next section consists in searching for suitable polynomials $P(x)$
and writing the corresponding rational period and continued fraction.

\section{Examples}

\begin{proposition} The condition of the proposition $P\in\Q[x]$ dividing
  $x^kP(x+1)+(x-1)^kP(x-1)$ is satisfied in the following cases:
  \begin{enumerate}
  \item $k\equiv0\pmod2$ and $P(x)=2x-1$.
  \item $k\equiv-1\pmod3$ and $P(x)=3x^2-3x+1$.
  \item $k\equiv-1\pmod4$ and $P(x)=2x^2-2x+1$.
  \item $k\equiv-1\pmod6$ and $P(x)=x^2-x+1$.
  \item $k=5$ and $P(x)=5x^4-10x^3+19x^2-14x+4$.
  \end{enumerate}
\end{proposition}

\Proof In the first four cases, it is sufficient to check that any root of
$P(x)=0$ is also a root of $x^kP(x+1)+(x-1)^kP(x-1)$. For instance in
the first case $P(x)=2x-1$, for $a=1/2$ we have $P(a+1)=2$ and $P(a-1)=-2$,
and indeed $2^{-k}\cdot 2+(-2)^{-k}\cdot(-2)=0$ when $k$ is even. The last
case is done by a direct divisibility test.\fp

The same proof shows the following:

\begin{proposition} The condition of the proposition $P\in\Q[x]$ dividing
  $x^kP(x+1)-(x-1)^kP(x-1)$ is satisfied in the following cases:
  \begin{enumerate}
  \item $k\equiv1\pmod2$ and $P(x)=2x-1$.
  \item $k\equiv1\pmod4$ and $P(x)=2x^2-2x+1$.
  \item $k\equiv2\pmod6$ and $P(x)=x^2-x+1$.
  \end{enumerate}
\end{proposition}

Note that I have not found any other examples than the ones given
in the above two propositions, but I may have missed some.

\smallskip

Thanks to these propositions, it is now just a matter of working out
explicitly all the above examples, in other words of computing the
partial fraction expansions of the expressions $1/(x^kP(x)P(x+1))$,
which is routine so not given explicitly. In particular, we will see that
in the non-alternating cases, the sum $S$ is a $\Q$-linear combination
of $1$ and $\z(k)$ for $k\ge2$ of fixed parity, and in the alternating
cases, with in addition $\log(2)$.

We now give the corresponding formulas, and give examples after.
The following corollary immediately follows from the above propositions:

\vfill\eject

\begin{corollary}
  By convention, set $\z(0)=\z^*(0)=0$, $\z^*(1)=\log(2)$, and
  $\z^*(k)=(2^{k-1}-1)\z(k)$ for $k\ge2$. We have the following general
  continued fractions:
\begin{align*}&\sum_{j=0}^{k-1}2^{2j}\z(2(k-j))=[[2^{2k-1},R_1(n)],[-1,-n^{4k}]]\;,\text{ with} \\   
  &R_1(x)=(x^{2k}(2x+1)+(x-1)^{2k}(2x-3))/(2x-1)\;,\\
&\sum_{j=0}^k(-3)^{3j}(\z(6(k-j)+2)+3\z(6(k-j)))=[[-(-3)^{3k+1}/2,R_2(n)],[1,-n^{12k+4}]]\;,\text{ with}\\
&R_2(x)=(x^{6k+2}(3x^2+3x+1)+(x-1)^{6k+2}(3x^2-9x+7))/(3x^2-3x+1)\;,\\
&\sum_{j=0}^k(-3)^{3j}(\z(6(k-j)+5)+3\z(6(k-j)+3))=[[(-3)^{3k+2}/2,R_3(n)],[1,-n^{12k+10}]]\;,\text{ with}\\
&R_3(x)=(x^{6k+5}(3x^2+3x+1)+(x-1)^{6k+5}(3x^2-9x+7))/(3x^2-3x+1)\;,\\
&\sum_{j=0}^k(-4)^j\z(4(k-j)+3)=[[4^k,R_4(n)],[1,-n^{8k+6}]]\;,\text{ with}\\
&R_4(x)=(x^{4k+3}(2x^2+2x+1)+(x-1)^{4k+3}(2x^2-6x+5))/(2x^2-2x+1)\;,\\
&\sum_{j=0}^k(\z(6(k-j)+5)-\z(6(k-j)+3))=[[-1/2,R_5(n)],[1,-n^{12k+10}]]\;,\text{ with}\\
&R_5(x)=(x^{6k+5}(x^2+x+1)+(x-1)^{6k+5}(x^3-3x+3))/(x^2-x+1)\;,\\
&\sum_{j=0}^k2^{4j}\z^*(2(k-j)+1)=[[2^{4k},R_6(n)],[-2^{2k},n^{4k+2}]]\;,\text{ with}\\
&R_6(x)=(x^{2k+1}(2x+1)-(x-1)^{2k+1}(2x-3))/(2x-1)\;,\\
&\sum_{j=0}^k(-64)^j\z^*(4(k-j)+1)=[[(-1)^k2^{6k-1},R_7(n)],[2^{4k},n^{8k+2}]]\;,\text{ with}\\
&R_7(x)=(x^{4k+1}(2x^2+2x+1)-(x-1)^{4k+1}(2x^2-6x+5))/(2x^2-2x+1)\;,\\
&\sum_{j=0}^k2^{6j}(\z^*(6(k-j)+2)-4\z^*(6(k-j)))=[[2^{6k},R_8(n)],[2^{6k+1},n^{12k+4}]]\;,\text{ with}\\
&R_8(x)=(x^{6k+2}(x^2+x+1)-(x-1)^{6k+2}(x^2-3x+3))/(x^2-x+1)\;.\end{align*}
\end{corollary}

We now give corresponding examples:

\begin{align*}
&\z(4)+4\z(2)=[[8,2n^4-4n^3+10n^2-8n+3],[-1,-n^8]]\;,\\
&27\z(2)-3\z(6)-\z(8)=[[81/2,R(n)],[-1,-n^{16}]]\;,\text{ with}\\
&R(n)=2n^8-8n^7+46n^6-110n^5+178n^4-182n^3+118n^2-44n+7\;,\\
&\z(5)+3\z(3)=[[9/2,2n^5-5n^4+22n^3-28n^2+23n-7],[1,-n^{10}]]\;,\\
&4\z(3)-\z(7)=[[4,2n^7-7n^6+37n^5-75n^4+99n^3-77n^2+31n-5],[-1,-n^{14}]]\;,\\
&\z(3)-\z(5)=[[1/2,2n^5-5n^4+22n^3-28n^2+15n-3],[-1,-n^{10}]]\;,\\
&4\z(5)+11\z(3)=[[273/16,2n^5-5n^4+42n^3-58n^2+45n-13],[4,-n^{10}]]\;,\\
&3\z(3)+16\log(2)=[[16,5n^2-5n+3],[-4,n^6]]\;,\\
&15\z(5)+48\z(3)+256\log(2)=[[256,7n^4-14n^3+18n^2-11n+3],[-16,n^{10}]]\;,\\
&64\log(2)-15\z(5)=[[32,9n^4-18n^3+30n^2-21n+5],[-16,n^{10}]]\;,\\
&127\z(8)-124\z(6)+64\z(2)=[[64,R(n)],[128,n^{16}]]\;,\text{ with}\\
&R(n)=12n^7-42n^6+110n^5-170n^4+154n^3-82n^2+24n-3\;.\end{align*}

\section{Second Method: use of the $\psi$ Function and Derivatives}

Recall that $\psi(z)$ is the logarithmic derivative of the gamma
function (most authors call $\psi$ the digamma function, and $\psi'$,
$\psi''$, etc... the trigamma, tetragamma functions, but this is terrible
terminology).

By orthogonality of characters, it is immediate to show that for $k\ge1$
$$\psi^{(k)}(r/m)=(-1)^{k-1}\dfrac{k!m^{k+1}}{\phi(m)}\sum_{\chi\bmod m}\ov{\chi}(r)L(\chi,k+1)\;,$$
and for $k=0$ the same formula is valid if we interpret $L(\chi_0,1)$ as
$$L(\chi_0,1)=\sum_{d\mid m}\dfrac{\mu(d)}{d}\log(d)-\dfrac{\phi(m)\log(m)}{m}\;.$$
In particular, for $m=1$, $2$, $3$, $4$, and $6$, which are the values of $m$
for which $\phi(m)\le2$, we obtain the following table, where as usual
we set $G=L(\chi_{-4},2)$, Catalan's constant, and $G_3=L(\chi_{-3},2)$:

\bigskip

\centerline{
\begin{tabular}{|c||c|c|c||}
  \hline
  $r/m$ & $\psi(r/m)+\ga$ & $\psi'(r/m)$ & $\psi''(r/m)$ \\
  \hline\hline
  1 & $0$ & $\z(2)$ & $-2\z(3)$ \\
  1/2 & $-2\log(2)$ & $3\z(2)$ & $-14\z(3)$ \\
  1/3 & $-3\log(3)/2-\pi/(2\sqrt{3})$ & $4\z(2)+9G_3/2$ & $-26\z(3)-4\pi^3/(3\sqrt{3})$ \\
  2/3 & $-3\log(3)/2+\pi/(2\sqrt{3})$ & $4\z(2)-9G_3/2$ & $-26\z(3)+4\pi^3/(3\sqrt{3})$ \\
  1/4 & $-3\log(2)-\pi/2$ & $6\z(2)+8G$ & $-56\z(3)-2\pi^3$ \\
  3/4 & $-3\log(2)+\pi/2$ & $6\z(2)-8G$ & $-56\z(3)+2\pi^3$ \\
  1/6 & $-\log(432)/2-3\pi/(2\sqrt{3})$ & $12\z(2)+45G_3/2$ & $-182\z(3)-12\pi^3/\sqrt{3}$ \\
  5/6 & $-\log(432)/2+3\pi/(2\sqrt{3})$ & $12\z(2)-45G_3/2$ & $-182\z(3)+12\pi^3/\sqrt{3}$ \\
  \hline
\end{tabular}}

\bigskip

On the other hand, there exist many continued fractions for $\psi(z)$ and
its derivatives. The ones for $\psi(z)$ itself are rather complicated,
and those for $\psi^{(k)}(z)$ for $k\ge3$ are trivial transformations of
the defining series, so the only remaining interesting ones are those for
$\psi'(z)$ and $\psi''(z)$. This of course implies that we restrict to
rational periods of degree two and three.

We choose the nicest continued fractions, taken from \cite{Cuyt}:

\medskip

For $\psi'(z)$ we have
$$\psi'(z)=[[0,(2z-1)(2n-1)],[2,n^4]]\;,$$
valid for $z>1/2$. However, from the trivial identity
$\psi'(z)=\psi'(z+1)+1/z^2$, we can
deduce infinitely many other continued fractions:
$$\psi'(z)=[[\sum_{0\le j<k}1/(z+j)^2,(2z+2k-1)(2n-1)],[2,n^4]]\;,$$
now valid for $z>1/2-k$.
Referring to the above table and choosing $z=1/3$, $2/3$, $1/4$, $3/4$,
$1/6$, and $5/6$, and $k=0$, $1$, etc..., we obtain as many continued
fractions as we like for $8\z(2)\pm9G_3$, $3\z(2)\pm4G$, and
$8\z(2)\pm15G_3$. For instance, after simplifications:
\begin{align*}
8\z(2)-9G_3&=[[0,2n-1],[12,9n^4]]\;,\\
8\z(2)+9G_3&=[[18,10n-5],[12,9n^4]]\;,\\
3\z(2)-4G&=[[0,2n-1],[2,4n^4]]\;,\\
3\z(2)+4G&=[[8,6n-3],[2,4n^4]]\;,\\
8\z(2)-15G_3&=[[0,4n-2],[4,9n^4]]\;,\\
8\z(2)+15G_3&=[[24,8n-4],[4,9n^4]]\;.\end{align*}

\medskip

For $\psi''(z)$ we have
$$\psi''(z)=[[0,(2n-1)(n^2-n+1+2z(z-1))],[-2,-n^6]]\;,$$
valid for $z>1/2$. However, from the trivial identity
$\psi''(z)=\psi''(z+1)-2/z^3$, we can
deduce infinitely many other continued fractions:
$$\psi''(z)=[[-2\sum_{0\le j<k}1/(z+j)^3,(2n-1)(n^2-n+1+2(z+k)(z+k-1))],[-2,-n^6]]\;,$$
now valid for $z>1/2-k$.
Referring to the above table and choosing $z=1/3$, $2/3$, $1/4$, $3/4$,
$1/6$, and $5/6$, and $k=0$, $1$, etc..., we obtain as many continued
fractions as we like for $39\z(3)\pm2\pi^3/\sqrt{3}$, $28\z(3)\pm\pi^3$,
and $91\z(3)\pm6\pi^3/\sqrt{3}$. For instance, after simplifications:
\begin{align*}
  39\z(3)-2\pi^3/\sqrt{3}&=[[0,(2n-1)(9n^2-9n+5)],[27,-81n^6]]\;,\\
  39\z(3)+2\pi^3/\sqrt{3}&=[[81,(2n-1)(9n^2-9n+17)],[27,-81n^6]]\;,\\
  28\z(3)-\pi^3&=[[0,(2n-1)(8n^2-8n+5)],[8,-64n^6]]\;,\\
  28\z(3)+\pi^3&=[[64,(2n-1)(8n^2-8n+13)],[8,-64n^6]]\;,\\
  91\z(3)-6\pi^3/\sqrt{3}&=[[0,(2n-1)(9n^2-9n+13/2)],[9,-81n^6]]\;,\\
  91\z(3)+6\pi^3/\sqrt{3}&=[[216,(2n-1)(9n^2-9n+25/2)],[9,-81n^6]]\;.
\end{align*}

\section{Third Method: Bauer--Muir Acceleration}

This very classical method is just as elementary as the previous ones, but
the formulas are slightly more complicated.

Let $(a(n),b(n))_{n\ge0}$ be a continued fraction with convergents
$(p(n),q(n))$, and let $r(n)_{n\ge1}$ be any sequence (for now). For
$n\ge1$ we define
$$R(n)=a(n)+r(n)\text{ and }d(n)=r(n)R(n+1)-b(n)=r(n)(a(n+1)+b(n+1))-b(n)\;.$$
We make the following two essential assumptions: $R(1)=a(1)+r(1)\ne0$, and
$d(n)\ne0$ for all $n\ge1$.

We define:
\begin{align*}
  A(0)&=a(0)+\dfrac{b(0)}{R(1)}\;,\quad B(0)=\dfrac{b(0)d(1)}{R(1)^2}\;,\quad
  A(1)=\dfrac{a(1)R(2)+b(1)}{R(1)}\;,\\
  A(n)&=R(n+1)-r(n-1)\dfrac{d(n)}{d(n-1)}\text{ for $n\ge2$, and }
  B(n)=b(n)\dfrac{d(n+1)}{d(n)}\text{ for $n\ge1$.}\end{align*}

The following result is easy to prove by induction:

\begin{proposition} Let $(P(n),Q(n))$ be the convergents of the continued
  fraction defined by $(A(n),B(n))$. For $n\ge2$ we have
  $$(P(n),Q(n))=(p(n+1),q(n+1))+r(n+1)(p(n),q(n))\;.$$
  In particular, if $p(n)/q(n)$ and $P(n)/Q(n)$ both tend to a limit as
  $n\to\infty$, these limits are equal.\end{proposition}

This process is called Bauer--Muir acceleration, because if $r(n)$ is chosen
appropriately, it accelerates the convergence of the continued fraction.
An important fact is that if the accelerated formulas are simple enough,
for instance when $d(n)$ is \emph{constant}, the acceleration process
can be \emph{iterated}. This fact, combined with a suitable diagonal process,
is the basis of Ap\'ery's initial proofs of the irrationality of $\z(2)$ and
$\z(3)$. However, we will not consider this here.

Let us consider some simple examples.

\bigskip

{\bf Example 1: $\log(2)$}

\smallskip

The trivial continued fraction for $\log(2)$, directly coming from the series
$\log(2)=\sum_{n\ge1}(-1)^{n-1}/n$, is $\log(2)=[[0,1],[1,n^2]]$.
Applying Bauer--Muir acceleration iteratively, we immediately obtain
$$\log(2)=[[0,1],[1,n^2]]=[[1,3],[-1,n^2]]=[[1/2,5],[1,n^2]]=[[5/6,7],[-1,n^2]]\;,$$
and so on, the general formula being
$$\log(2)=[[\sum_{1\le j\le k}(-1)^{j-1}/j,2k+1],[(-1)^k,n^2]]\;.$$
Note that this is \emph{not} the same as the trivial continued fraction
obtained from $\sum_{n>k}(-1)^{n-1}/n$, since this converges like $(-1)^n/n$
as the initial one, while the accelerated formula given above converges
like $(-1)^n/n^{2k+1}$.

\bigskip

{\bf Example 2: $\pi^2/6$}

\smallskip

The trivial continued fraction for $\z(2)=\pi^2/6$, directly coming from the
series $\pi^2/6=\sum_{n\ge1}1/n^2$, is $\pi^2/6=[[0,2n^2-2n+1],[1,-n^4]]$.
Applying Bauer--Muir acceleration iteratively, we immediately obtain
\begin{align*}\pi^2/6
  &=[[0,2n^2-2n+1],[1,-n^4]]=[[2,2n^2-2n+3],[-1,-n^4]]\\
  &=[[3/2,2n^2-2n+7],[1,-n^4]]=[[31/18,2n^2-2n+13],[-1,-n^4]]\\
  &=[[115/72,2n^2-2n+21],[1,-n^4]]=[[3019/1800,2n^2-2n+31],[-1,-n^4]]\\
  &=[[973/600,2n^2-2n+43],[1,-n^4]]\;,
\end{align*}
and so on, the general formula being
$$\pi^2/6=[[2\sum_{1\le j\le k}(-1)^{j-1}/j^2,2n^2-2n+k^2+k+1],[(-1)^k,-n^4]]\;.$$
Once again, this is not the tail of the series defining $\pi^2/6$, since
it converges like $1/n^{2k+1}$.

\bigskip

{\bf Example 3: $\pi^2/6$ (again)}

\smallskip

Another trivial continued fraction for $\pi^2/6$, directly coming from the
series $\pi^2/6=2\sum_{n\ge1}(-1)^{n-1}/n^2$, is
$\pi^2/6=[[0,2n-1],[2,n^4]]$.
Applying Bauer--Muir acceleration iteratively, we immediately obtain
\begin{align*}\pi^2/6
  &=[[0,2n-1],[2,n^4]]=[[1,6n-3],[2,n^4]]=[[5/4,10n-5],[2,n^4]]\\
  &=[[49/36,14n-7],[2,n^4]]=[[205/144,18n-9],[2,n^4]]\\
  &=[[5269/3600,22n-11],[2,n^4]]\;,\end{align*}
and so on, the general formula being
$$\pi^2/6=[[\sum_{1\le j\le k}1/j^2,(2k+1)(2n-1)],[2,n^4]]\;.$$
Note that the constant term of the continued fraction for $\sum_{j\ge1}1/j^2$
is the partial sum of the series $\sum_{j\ge1}(-1)^{j-1}/j^2$ and the
constant term of the continued fraction for $\sum_{j\ge1}(-1)^{j-1}/j^2$ is
the partial sum of the series $\sum_{j\ge1}1/j^2$.

\bigskip

{\bf Example 4: $G=L(\chi_{-4},2)$, Catalan's constant}

\smallskip
Here, we could take the trivial continued fraction for $G$, directly coming
from the series $G=\sum_{n\ge1}(-1)^{n-1}/(2n-1)^2$, and apply iteratively
Bauer--Muir, giving
\begin{align*}G
  &=[[0,1,8(n-1)],[1,(2n-1)^4]]=[[1/6,7/3,24(n-1)],[16/9,(2n-1)^4]]\\
  &=[[19/82,145/41,40(n-1)],[4096/1681,(2n-1)^4]]\;,\end{align*}
and so on, but this is not pretty, first because $b(0)$ changes, and second
because one needs to specify both $a(0)$ and $a(1)$ since $n-1$ vanishes
for $n=1$.

\smallskip

A nicer continued fraction taken from \cite{Cuyt}, which in fact we are
going to prove, is $G=[[1,8n^2-8n+7],[-1/2,-16n^4]]$. We thus obtain
\begin{align*}G&=[[0,8n^2-8n+3],[1/2,-16n^4]]\\
  &=[[1,8n^2-8n+7],[-1/2,-16n^4]]=[[8/9,8n^2-8n+19],[1/2,-16n^4]]\\
  &=[[209/225,8n^2-8n+39],[-1/2,-16n^4]]\\
  &=[[10016/11025,8n^2-8n+67],[1/2,-16n^4]]\;,\end{align*}
and so on, the general formula being
$$G=[[\sum_{1\le j\le k}(-1)^{j-1}/(2j-1)^2,8n^2-8n+4k^2+3],[(-1)^k/2,-16n^4]]\;,$$
and the continued fraction for $k=0$ being obtained by \emph{reverse}
Bauer--Muir, with extremely slow convergence in $1/\log(n)$.

But in turn these formulas \emph{prove} that the initial continued fraction
converges to $G$, since if we set $S_k=\sum_{1\le j\le k}(-1)^{j-1}/(2j-1)^2$,
which is the $k$th partial sum of the series defining $G$, the $k$th continued
fraction is $S_k+(-1)^k/2/(4k^2+3-16/(4k^2+19-\cdots))$, and this clearly tends
to $\lim_{k\to\infty}S_k=G$ as $k\to\infty$.

\bigskip

{\bf Example 5: $\z(3)$}

\smallskip

The trivial continued fraction for $\z(3)$, directly coming from the
series $\z(3)=\sum_{n\ge1}1/n^3$, is
$\z(3)=[[0,(2n-1)(n^2-n+1)],[1,-n^6]]$.
Applying Bauer--Muir acceleration iteratively, we immediately obtain
\begin{align*}\z(3)
  &=[[0,(2n-1)(n^2-n+1)],[1,-n^6]]=[[1,(2n-1)(n^2-n+5)],[1,-n^6]]\\
  &=[[9/8,(2n-1)(n^2-n+13)],[1,-n^6]]\\
  &=[[251/216,(2n-1)(n^2-n+25)],[1,-n^6]]\;,
\end{align*}
and so on, the general formula being
$$\z(3)=[[\sum_{1\le j\le k}1/j^3,(2n-1)(n^2-n+2k^2+2k+1)],[1,-n^6]]\;.$$

\smallskip

The reader can check that unfortunately, the method does not work (i.e.,
the formulas become extremely complicated) for the
alternating sum giving $3\z(3)/4$, nor for $\z(k)$ for $k\ge4$, explaining
in large part why Ap\'ery's method has not been extended to $\z(k)$ for
$k\ge5$. Note that, on the contrary, Ap\'ery's method does work for
a large number of other series, and will be the object of a future
paper.

\section{Conclusion}

We have given three rather different methods to obtain infinitely many
continued fractions for certain linear combinations of zeta and $L$ values.
Note, however, that they are all polynomially convergent
(i.e., in $C/n^k$ for some $k\ge1$), while really interesting continued
fractions are exponentially or at least sub-exponentially convergent.
As already mentioned, this will be the subject of a future paper \cite{Coh}.

\end{document}